\documentclass{amsart}

\usepackage{fancyhdr}
\usepackage{lastpage}
\usepackage{stmaryrd,yhmath}

\newcommand{\bdis}{\begin{displaymath}}
\newcommand{\edis}{\end{displaymath}}
\newcommand{\be}{\begin{equation}}
\newcommand{\ee}{\end{equation}}
\newcommand{\mbb}{\mathbb}
\newcommand{\mcal}{\mathcal}

\newcommand{\vp}{\varphi}

\newcommand{\mG}{\mathring{G}}
\newcommand{\mT}{\mathring{T}}
\newcommand{\mg}{\mathring{g}}

\newcommand{\zf}{\zeta\left(\frac{1}{2}+it\right)}

\theoremstyle{definition}

\newtheorem{cor}[]{Corollary}

\theoremstyle{remark}
\newtheorem{remark}[]{Remark}

\newtheorem*{mydef1}{{\bf Theorem}}

\numberwithin{equation}{section}

\begin{document}

\title{Jacob's ladders and some new consequences from A. Selberg's formula}

\author{Jan Moser}

\address{Department of Mathematical Analysis and Numerical Mathematics, Comenius University, Mlynska Dolina M105, 842 48 Bratislava, SLOVAKIA}

\email{jan.mozer@fmph.uniba.sk}

\keywords{Riemann zeta-function}

\begin{abstract}
It is proved in this paper that the Jacob's ladders together with the A. Selberg's classical formula (1942) lead to a new kind of formulae for some
short trigonometric sums. These formulae cannot be obtained in the classical theory of A. Selberg, and all the less, in the theories of Balasubramanian,
Heath-Brown and Ivic.
\end{abstract}

\maketitle

\section{The A. Selberg's formula}

A. Selberg has proved in 1942 the following formula
\be \label{1.1}
\int_T^{T+U}X^2(t)\left(\frac{n_2}{n_1}\right)^{it}{\rm d}t=\sqrt{\frac \pi 2}\frac{U}{\sqrt{n_1n_2}}\left(\ln\frac{P^2}{n_1n_2}+2c\right)+
\mcal{O}(T^{1/2}\xi^5)
\ee
(see \cite{19}, p. 55), where
\be \label{1.2}
\begin{split}
& X(t)=\frac 12 t^{1/4}e^{\frac{1}{4}\pi t}\pi^{-\frac s2}\zeta(s),\ s=\frac 12 + it , \\
& U=T^{1/2+\epsilon},\ \xi=\left(\frac{T}{2\pi}\right)^{\epsilon/10},\ \epsilon\leq \frac{1}{10},\ P=\sqrt{\frac{T}{2\pi}}\\
& n_1,n_2\in \mbb{N}, (n_1,n_2)=1,\ n_1,n_2\leq\xi ,
\end{split}
\ee
(comp. \cite{19}, pp. 10, 18, $a=1/2+\epsilon,\ \epsilon>0$) and $c$ is the Euler's constant. Since (see \cite{19}, p. 10, \cite{20}, p. 79)
\bdis
Z^2(t)=\left|\zf\right|^2=\sqrt{\frac 2\pi}X^2(t)\left( 1+\mcal{O}(\frac 1t)\right),
\edis
i.e.
\be \label{1.3}
X^2(t)=\sqrt{\frac 2\pi}Z^2(t)\left( 1+\mcal{O}(\frac 1t)\right)
\ee
where
\be \label{1.4}
\begin{split}
& Z(t)=e^{i\vartheta(t)}\zf , \\
& \vartheta(t)=-\frac 12t\ln\pi+\text{Im}\ln\Gamma\left(\frac 14+\frac 12 it\right)=\frac t2\ln\frac{t}{2\pi}-\frac t2-\frac \pi 8 +
\mcal{O}(\frac 1t)
\end{split}
\ee
is the signal defined by the Riemann zeta-function $\zeta(s)$. Following eqs. (\ref{1.1}) and (\ref{1.3}) we obtain
\be \label{1.5}
\int_T^{T+U} Z^2(t)\left(\frac{n_2}{n_1}\right)^{it}{\rm d}t=\frac{U}{\sqrt{n_1n_2}}\left(\ln\frac{P^2}{n_1n_2}+2c\right)+
\mcal{O}(T^{1/2}\xi^5)
\ee

\begin{remark}
If $n_1=n_2=1$ then the Hardy-Littlewood-Ingham formula
\bdis
\int_T^{T+U} Z^2(t){\rm d}t=U\ln\frac{T}{2\pi}+2cU+\mcal{O}(T^{1/2}\xi^5)
\edis
follows from the A. Selberg's formula (\ref{1.5}) (comp. \cite{20}, p. 120).
\end{remark}

\begin{remark}
Let us remind that the A. Selberg's formula (\ref{1.5}) played the main role in proving the fundamental Selberg's result
\bdis
N_0(T+U)-N_0(T)>A(\epsilon)U\ln T
\edis
where $N_0$ stands for the number of zeroes of the function $\zeta(1/2+it),\ t\in (0,T]$.
\end{remark}

In this paper it is proved that the Jacob's ladders together with the A. Selberg's classical formula lead to a new kind of results for some short
trigonometric sums.

This paper is a continuation of the series of works \cite{3} - \cite{18}.

\section{The result}

\subsection{}

Let us remind some notions. First of all
\be \label{2.1}
\tilde{Z}^2(t)=\frac{{\rm d}\vp_1(t)}{{\rm d}t},\ \vp_1(t)=\frac 12\vp(t) ,
\ee
where
\be \label{2.2}
\tilde{Z}^2(t)=\frac{Z^2(t)}{2\Phi^\prime_\vp[\vp(t)]}=\frac{Z^2(t)}{\left\{ 1+\mcal{O}\left(\frac{\ln\ln t}{\ln t}\right)\right\}\ln t}
\ee
(see \cite{3}, (3.9); \cite{5}, (1.3); \cite{9}, (1.1), (3.1), (3.2)) and $\vp(t)$ is the Jacob's ladder, i.e. the solution of the following nonlinear
integral equation
\bdis
\int_0^{\mu[x(T)]}Z^2(t)e^{-\frac{2}{x(T)}t}{\rm d}t=\int_0^T Z^2(t){\rm d}t
\edis
that was introduced in our paper \cite{3}. Next, we have (see \cite{1}, comp. \cite{18})
\be \label{2.3}
\begin{split}
& G_3(x)=G_3(x;T,U)= \\
& =\bigcup_{T\leq g_{2\nu}\leq T+U}\{ t:\ g_{2\nu}(-x)\leq t\leq g_{2\nu}(x)\},\ 0<x\leq \frac \pi 2 , \\
& G_4(y)=G_4(y;T,U)= \\
& = \bigcup_{T\leq g_{2\nu+1}\leq T+U}\{ t:\ g_{2\nu+1}(-y)\leq t\leq g_{2\nu+1}(y)\},\ 0<y\leq \frac \pi 2 ,
\end{split}
\ee
and the collection of sequences $\{ g_\nu(\tau)\},\ \tau\in [-\pi,\pi],\ \nu=1,2,\dots $ is defined by the equation (see \cite{1}, \cite{18}, (6))
\bdis
\vartheta_1[g_\nu(\tau)]=\frac \pi 2\nu+\frac{\tau}{2};\ g_\nu(0)=g_\nu
\edis
where (comp. (\ref{1.4}))
\bdis
\vartheta_1(t)=\frac t2\ln\frac{t}{2\pi}-\frac t2-\frac \pi 8 .
\edis

\subsection{}

In this paper we obtain some new integrals containing the following short trigonometric sums
\bdis
\begin{split}
& \sum_{2\leq p\leq\xi}\frac{1}{\sqrt{p}}\cos(t\ln p),\ \sum_{2\leq n\leq\xi}\frac{1}{\sqrt{n}}\cos(t\ln n) , \\
& \sum_{2\leq n\leq\xi}\frac{d(n)}{\sqrt{n}}\cos(t\ln n)
\end{split}
\edis
where $p$ is the prime, $n\in\mbb{N}$ and $d(n)$ is the number of divisors of $n$. In this direction, the following theorem holds true.

\begin{mydef1}
Let
\be \label{2.4}
G_3(x)=\vp_1(\mathring{G}_3(x)),\ G_4(y)=\vp_1(\mathring{G}_4(y)) .
\ee
Then we have
\be \label{2.5}
\begin{split}
& \int_{\mG_3(x)}\left(\sum_{2\leq p\leq\xi}\frac{1}{\sqrt{p}}\cos(\vp_1(t)\ln p)\right)Z^2\{\vp_1(t)\}\tilde{Z}^2(t){\rm d}t\sim \\
& \frac{2x}{\pi}U\ln P\ln\ln P,\ x\in (0,\pi/2] , \\
& \int_{\mG_4(y)}\left(\sum_{2\leq p\leq\xi}\frac{1}{\sqrt{p}}\cos(\vp_1(t)\ln p)\right)Z^2\{\vp_1(t)\}\tilde{Z}^2(t){\rm d}t\sim \\
& \frac{2y}{\pi}U\ln P\ln\ln P,\ y\in (0,\pi/2] ,
\end{split}
\ee

\be \label{2.6}
\begin{split}
& \int_{\mG_3(x)}\left(\sum_{2\leq n\leq\xi}\frac{1}{\sqrt{n}}\cos(\vp_1(t)\ln n)\right)Z^2\{\vp_1(t)\}\tilde{Z}^2(t){\rm d}t\sim \\
& \frac{1}{\pi}\left\{ \left(\frac{2\epsilon}{5}-\frac{\epsilon^2}{50}\right)x+\frac{\epsilon^2}{50}\sin x\right\}U\ln^2 P , \\
& \int_{\mG_4(y)}\left(\sum_{2\leq n\leq\xi}\frac{1}{\sqrt{n}}\cos(\vp_1(t)\ln n)\right)Z^2\{\vp_1(t)\}\tilde{Z}^2(t){\rm d}t\sim \\
& \frac{1}{\pi}\left\{ \left(\frac{2\epsilon}{5}-\frac{\epsilon^2}{50}\right)y-\frac{\epsilon^2}{50}\sin y\right\}U\ln^2 P ,
\end{split}
\ee

\be \label{2.7}
\begin{split}
& \int_{\mG_3(x)}\left(\sum_{2\leq n\leq\xi}\frac{d(n)}{\sqrt{n}}\cos(\vp_1(t)\ln n)\right)Z^2\{\vp_1(t)\}\tilde{Z}^2(t){\rm d}t\sim \\
& \frac{\sin x}{2500\pi^3}U\ln^4P , \\
& \int_{\mG_4(y)}\left(\sum_{2\leq n\leq\xi}\frac{d(n)}{\sqrt{n}}\cos(\vp_1(t)\ln n)\right)Z^2\{\vp_1(t)\}\tilde{Z}^2(t){\rm d}t\sim \\
& -\frac{\sin y}{2500\pi^3}U\ln^2P ,
\end{split}
\ee
where
\be \label{2.8}
t-\vp_1(t) \sim (1-c)\pi(t),\ t\to\infty ,
\ee
and $\pi(t)$ is the prime-counting function.
\end{mydef1}

\begin{remark}
Let $T=\vp_1(\mT)$, $T+U=\vp_1(\widering{T+U})$, (comp. (\ref{2.4})). Then from (\ref{2.8}), similarly to \cite{14}, (1.8), we obtain
\bdis
\rho\{ [T,T+U]; [\mT,\widering{T+U}]\}\sim (1-c)\pi(T);\ T+U<\mT,
\edis
where $\rho$ stands for the distance of the corresponding segments.
\end{remark}

\begin{remark}
The formulae (\ref{2.5}) - (\ref{2.7}) cannot be obtained in the classical theory of A. Selberg, and, all the less, in the theories of Balasubramanian,
Heath-Brown and Ivic.
\end{remark}

\section{New asymptotic formulae for the short trigonometric sums: their dependence on $|\zf|^2$}

We obtain, putting $x=y=\pi/2$ in (\ref{2.5})
\bdis
\begin{split}
& \int_{\mG_3(\pi/2)\cup\mG_4(\pi/2)}\left(\sum_{2\leq p\leq \xi}\frac{1}{\sqrt{p}}\cos\{\vp_1(t)\ln p\}\right)
Z^2\{\vp_1(t)\}\tilde{Z}^2(t){\rm d}t\sim \\
& 2U\ln P\ln\ln P .
\end{split}
\edis
Using successively the mean-value theorem (since $\mG_3(\pi/2)\cup\mG_4(\pi/2)$ is a segment), we have
\be \label{3.1}
\begin{split}
& \sum_{2\leq p\leq \xi}\frac{1}{\sqrt{p}}\cos\{\vp_1(\alpha_1)\ln p\}\int_{\mG_3(\pi/2)\cup\mG_4(\pi/2)}Z^2\{\vp_1(t)\}Z^2(t){\rm d}t= \\
& =\sum_{2\leq p\leq \xi}\frac{1}{\sqrt{p}}\cos\{\vp_1(\alpha_1)\ln p\}Z^2\{\vp_1(\alpha_2)\}
\int_{\mG_3(\pi/2)\cup\mG_4(\pi/2)}\tilde{Z}^2(t){\rm d}t\sim \\
& 2U\ln P\ln\ln P,\ \alpha_1,\alpha_2\in \mG_3(\pi/2)\cup\mG_4(\pi/2);\ \alpha_1=\alpha_1(T,U)=\alpha_1(T,\epsilon), \dots  .
\end{split}
\ee
Since
\be \label{3.2}
\int_{\mG_3(\pi/2)\cup\mG_4(\pi/2)}\tilde{Z}^2{\rm d}t=\left|\mG_3(\pi/2)\cup\mG_4(\pi/2)\right|
\ee
(comp. Remark 8), and
\be \label{3.3}
m\{\mG_3(x)\}\sim \frac x\pi U,\ m\{\mG_4(y)\}\sim \frac y\pi U \ \Rightarrow \ \left|\mG_3(\pi/2)\cup\mG_4(\pi/2)\right|\sim U
\ee
(see \cite{2}, (13), $m$ stands for the measure) then we obtain from (\ref{2.5}) (see (\ref{3.1}) - (\ref{3.3})) the following

\begin{cor}
For every $T\geq T_0[\vp_1]$ there are the values $\alpha_1(T), \alpha_2(T)\in \mG_3(\pi/2)\cup\mG_4(\pi/2)$ such that
\be \label{3.4}
\sum_{2\leq p\leq \xi}\frac{1}{\sqrt{p}}\cos\{\vp_1(\alpha_1(T))\ln p\}\sim \frac{2\ln P\ln\ln P}{\left|\zeta\left(\frac{1}{2}+i
\vp_1(\alpha_2(T))\right)\right|^2},\ T\to\infty
\ee
where $\vp_1(\alpha_1(T)), \vp_1(\alpha_2(T))\in \mG_3(\pi/2)\cup\mG_4(\pi/2)$.
\end{cor}

Similarly, we obtain from (\ref{2.6})

\begin{cor}
For every $T\geq T_0[\vp_1]$ there are the values $\alpha_3(T), \alpha_4(T)\in \mG_3(\pi/2)\cup\mG_4(\pi/2)$ such that
\be \label{3.5}
\begin{split}
& \sum_{2\leq n\leq \xi}\frac{1}{\sqrt{n}}\cos\{\vp_1(\alpha_3(T))\ln n\}\sim \\
& \sim \left(\frac{2\epsilon}{5}-\frac{\epsilon^2}{50}\right)\frac{\ln^2 P}
{\left|\zeta\left(\frac 12 +i\vp_1(\alpha_4(T))\right)\right|^2},\ T\to\infty
\end{split}
\ee
where $\vp_1(\alpha_3(T)), \vp_1(\alpha_4(T))\in G_3(\pi/2)\cup G_4(\pi/2)$.
\end{cor}

\begin{remark}
From the asymptotic formulae (\ref{3.4}), (\ref{3.5}) it follows that the values of mentioned short trigonometric sums are connected with the values
of the Riemann zeta-function $\zf$ for some infinite subset of $t$.
\end{remark}

\section{New asymptotic formulae on two collections of disconnected sets $G_3(x), G_4(y)$}

From (\ref{2.7}), similarly to p. 3, we obtain

\begin{cor}
\be \label{4.1}
\begin{split}
& \left.\left\langle\sum_{2\leq n\leq \xi}\frac{d(n)}{\sqrt{n}}\cos\{\vp_1(t)\ln n\}\right\rangle\right|_{\mG_3(x)}\sim \\
& \sim \frac{1}{2500\pi^2}\frac{\sin x}{x}\frac{\ln^4P}{\langle Z^2\{\vp_1(t)\}\rangle|_{\mG_3(x)}} \\
& \left.\left\langle\sum_{2\leq n\leq \xi}\frac{d(n)}{\sqrt{n}}\cos\{\vp_1(t)\ln n\}\right\rangle\right|_{\mG_4(y)}\sim \\
& \sim -\frac{1}{2500\pi^2}\frac{\sin y}{y}\frac{\ln^4P}{\langle Z^2\{\vp_1(t)\}\rangle|_{\mG_4(y)}},\ T\to\infty
\end{split}
\ee
where $\langle(\dots )\rangle|_{\mG_3(x)},\dots$ denote the mean-value of $(\dots)$ on $\mG_3(x), \dots$ .
\end{cor}

\begin{remark}
It follows from (\ref{4.1}) that the short trigonometric sum
\bdis
\sum_{2\leq n\leq \xi}\frac{d(n)}{\sqrt{n}}\cos\{t\ln n\},\ t\geq T_0[\vp_1]
\edis
has an infinitely many zeroes of the odd order.
\end{remark}

\section{Law of the asymptotic equality of areas}

Let
\bdis
\begin{split}
& \mG_3^+(x)=\left\{ t:\ t\in \mG_3(x),\ \sum_{2\leq n\leq \xi}\frac{d(n)}{\sqrt{n}}\cos\{\vp_1(t)\ln n\}>0\right\}, \\
& \vdots \\
& \mG_4^-(x)=\left\{ t:\ t\in \mG_4(x),\ \sum_{2\leq n\leq \xi}\frac{d(n)}{\sqrt{n}}\cos\{\vp_1(t)\ln n\}<0\right\} .
\end{split}
\edis
Then we obtain from (\ref{2.7}), (comp. Corollary 3 in \cite{14})

\begin{cor}
\be \label{5.1}
\begin{split}
& \int_{\mG_3^+(x)\cup \mG_4^+(x)}\left(\sum_{2\leq n\leq \xi}\frac{d(n)}{\sqrt{n}}\cos\{\vp_1(t)\ln n\}\right)
Z^2\{\vp_1(t)\}\tilde{Z}^2(t){\rm d}t\sim \\
& \sim -\int_{\mG_3^-(x)\cup \mG_4^-(x)}\left(\sum_{2\leq n\leq \xi}\frac{d(n)}{\sqrt{n}}\cos\{\vp_1(t)\ln n\}\right)
Z^2\{\vp_1(t)\}\tilde{Z}^2(t){\rm d}t .
\end{split}
\ee
\end{cor}

\begin{remark}
The formula (\ref{5.1}) represents the law of the asymptotic equality of the areas (measures) of complicated figures corresponding to the positive part and
the negative part, respectively, of the graph of the function
\be \label{5.2}
\sum_{2\leq n\leq \xi}\frac{d(n)}{\sqrt{n}}\cos\{\vp_1(t)\ln n\}Z^2\{\vp_1(t)\}\tilde{Z}^2(t),\ t\in\mG_3(x)\cup\mG_4(x) ,
\ee
where $x\in (0,\pi/2]$. This is one of the laws governing the \emph{chaotic} behaviour of the positive and negative values of the signal (\ref{5.2}).
This signal is created by the complicated modulation of the fundamental signal $Z(t)=e^{i\vartheta(t)}\zf$, (comp. (\ref{1.4}), (\ref{2.2})).
\end{remark}

\section{Proof of the Theorem}

\subsection{}

Let us remind that the following lemma holds true (see \cite{8}, (2.5); \cite{9}, (3.3)): for every integrable function (in the Lebesgue sense)
$f(x), \ x\in [\vp_1(T),\vp_1(T+U)]$ we have
\be \label{6.1}
\int_T^{T+U} f[\vp_1(t)]\tilde{Z}^2(t){\rm d}t=\int_{\vp_1(T)}^{\vp_1(T+U)}f(x){\rm d}x,\ U\in (0,T/\ln T] ,
\ee
where $t-\vp_1(t)\sim (1-c)\pi(t)$. In the case (comp. (\ref{2.4})) $T=\vp_1(\mT)$, $T+U=\vp_1(\widering{T+U})$, we obtain from (\ref{6.1}) the following
equality
\be \label{6.2}
\int_{\mT}^{\widering{T+U}}f[\vp_1(t)]\tilde{Z}^2(t){\rm d}t=\int_T^{T+U}f(x){\rm d}x .
\ee

\subsection{}

First of all, we have from (\ref{6.2}), for example,
\bdis
\int_{\mg_{2\nu}(-x)}^{\mg_{2\nu}(x)}f[\vp_1(t)]\tilde{Z}^2(t){\rm d}t=\int_{g_{2\nu}(-x)}^{g_{2\nu}(x)}f(t){\rm d}t ,
\edis
(see (\ref{2.3}). Next, in the case
\bdis
f(t)=\left(\sum_{2\leq p\leq \xi}\frac{1}{\sqrt{p}}\cos\{\vp_1(t)\ln p\}\right)Z^2\{\vp_1(t)\}
\edis
we have the following $\tilde{Z}^2$-transformation
\be \label{6.3}
\begin{split}
& \int_{\mG_3(x)}\left(\sum_{2\leq p\leq \xi}\frac{1}{\sqrt{p}}\cos\{\vp_1(t)\ln p\}\right)Z^2\{\vp_1(t)\}\tilde{Z}^2(t){\rm d}t= \\
& = \int_{\mG_3(x)}\left(\sum_{2\leq p\leq \xi}\frac{1}{\sqrt{p}}\cos\{\vp_1(t)\ln p\}\right)Z^2(t){\rm d}t , \\
& \int_{\mG_4(y)}\left(\sum_{2\leq p\leq \xi}\frac{1}{\sqrt{p}}\cos\{\vp_1(t)\ln p\}\right)Z^2\{\vp_1(t)\}\tilde{Z}^2(t){\rm d}t= \\
& = \int_{\mG_4(y)}\left(\sum_{2\leq p\leq \xi}\frac{1}{\sqrt{p}}\cos\{\vp_1(t)\ln p\}\right)Z^2(t){\rm d}t .
\end{split}
\ee
Let us remind that we have proved (see \cite{2}, (13) and Corollary 7) the following formulae
\be \label{6.4}
\begin{split}
& \int_{G_3(x)}\left(\sum_{2\leq p\leq \xi}\frac{1}{\sqrt{p}}\cos\{\vp_1(t)\ln p\}\right)Z^2(t){\rm d}t\sim \frac{2x}{\pi}\ln P\ln\ln P, \\
& \int_{G_4(y)}\left(\sum_{2\leq p\leq \xi}\frac{1}{\sqrt{p}}\cos\{\vp_1(t)\ln p\}\right)Z^2(t){\rm d}t\sim \frac{2y}{\pi}\ln P\ln\ln P .
\end{split}
\ee
Now, our formulae (\ref{2.5}) follow from (\ref{6.3}), (\ref{6.4}).

\subsection{}

Similarly, from the formulae
\bdis
\begin{split}
& \int_{G_3(x)}\left(\sum_{2\leq n\leq \xi}\frac{1}{\sqrt{n}}\cos\{t\ln n\}\right)Z^2(t){\rm d}t \sim \\
& \sim \frac x\pi\left(\frac{2\epsilon}{5}-\frac{\epsilon^2}{50}+\frac{\epsilon^2}{50}\frac{\sin x}{x}\right)U\ln^2P , \\
& \int_{G_4(y)}\left(\sum_{2\leq n\leq \xi}\frac{1}{\sqrt{n}}\cos\{t\ln n\}\right)Z^2(t){\rm d}t \sim \\
& \sim \frac y\pi\left(\frac{2\epsilon}{5}-\frac{\epsilon^2}{50}-\frac{\epsilon^2}{50}\frac{\sin y}{y}\right)U\ln^2P ,
\end{split}
\edis
and
\bdis
\begin{split}
& \int_{G_3(x)}\left(\sum_{2\leq n\leq \xi}\frac{d(n)}{\sqrt{n}}\cos\{t\ln n\}\right)Z^2(t){\rm d}t \sim \frac{\sin x}{2500\pi^3}U\ln^4P , \\
& \int_{G_4(y)}\left(\sum_{2\leq n\leq \xi}\frac{d(n)}{\sqrt{n}}\cos\{t\ln n\}\right)Z^2(t){\rm d}t \sim -\frac{\sin y}{2500\pi^3}U\ln^4P
\end{split}
\edis
(see \cite{2}, (13) and Corollaries 8 and 9) we obtain (\ref{2.6}) and (\ref{2.7}), respectively.

\begin{remark}
The formulae of type (\ref{3.2}) can be obtained from (\ref{6.2}) putting $f(t)\equiv 1$.
\end{remark}

\thanks{I would like to thank Michal Demetrian for helping me with the electronic version of this work.}


\begin{thebibliography}{29}
\bibitem{18}
J. Moser, `Jacob's ladders and the nonlocal interaction of the function $Z^2(t)$ with the function $\tilde{Z}^2(t)$ on the distance $\sim (1-c)\pi(t)$
for the collections of disconnected sets`, (2010), arXiv: 1007.5147.
\end{thebibliography}
\end{document}